\documentclass[11pt, 
 reprint,
 notitlepage,
 amsmath,amssymb,
 aps,pra,superscriptaddress,
]{revtex4-2} 

\usepackage[T1]{fontenc}
\usepackage[colorlinks,citecolor=blue,urlcolor=blue]{hyperref} 
\usepackage{graphics}
\usepackage{natbib}
\usepackage{physics}
\usepackage{siunitx}
\usepackage{listings} 
\usepackage{xcolor}

\definecolor{codegreen}{rgb}{0,0.6,0}
\definecolor{codegray}{rgb}{0.5,0.5,0.5}
\definecolor{codepurple}{rgb}{0.58,0,0.82}
\definecolor{backcolour}{rgb}{0.95,0.95,0.92}

\lstdefinestyle{mystyle}{
    backgroundcolor=\color{backcolour},   
    commentstyle=\color{codegreen},
    keywordstyle=\color{magenta},
    numberstyle=\tiny\color{codegray},
    stringstyle=\color{codepurple},
    basicstyle=\ttfamily\footnotesize,
    breakatwhitespace=false,         
    breaklines=true,                 
    captionpos=b,                    
    keepspaces=true,                 
    numbers=left,                    
    numbersep=5pt,                  
    showspaces=false,                
    showstringspaces=false,
    showtabs=false,                  
    tabsize=2
}

\begin{document}

\title{Exponential Integrators for Phase-Field Equations using Pseudo-spectral Methods: A Python Implementation}
\author{Elvis do A. Soares}%
\email{elvis.asoares@gmail.com}
\author{Amaro G. Barreto Jr.}%
\affiliation{Engenharia de Processos Químicos e Bioquímicos (EPQB), Escola de Química, Universidade Federal do Rio de Janeiro, 21941-909, Rio de Janeiro, RJ, Brazil}%
\author{Frederico W. Tavares}%
\email{tavares@eq.ufrj.br}
\affiliation{Engenharia de Processos Químicos e Bioquímicos (EPQB), Escola de Química, Universidade Federal do Rio de Janeiro, 21941-909, Rio de Janeiro, RJ, Brazil}%
\affiliation{Programa de Engenharia Química, COPPE, Universidade Federal do Rio de Janeiro, 21941-909, Rio de Janeiro, RJ, Brazil}%

\date{\today}

\begin{abstract}
    In this paper, we implement exponential integrators, specifically Integrating Factor (IF) and Exponential Time Differencing (ETD) methods, using pseudo-spectral techniques to solve phase-field equations within a Python framework. These exponential integrators have showcased robust performance and accuracy when addressing stiff nonlinear partial differential equations. We compare these integrators to the well-known implicit-explicit (IMEX) Euler integrators used in phase-field modeling. The synergy between pseudo-spectral techniques and exponential integrators yields significant benefits for modeling intricate systems governed by phase-field dynamics, such as solidification processes and pattern formation. Our comprehensive Python implementation illustrates the effectiveness of this combined approach in solving phase-field model equations. The results obtained from this implementation highlight the accuracy and computational advantages of the ETD method compared to other numerical techniques.
\end{abstract}

\maketitle

\section{Introduction}

Phase-field equations~\cite{Chen2002,Gallego2011,Chen2022} are popular approaches for understanding and modeling complex dynamics, such as phase separation and solidification processes. The Cahn-Hilliard (CH) model~\cite{Cahn1958,Cahn1959,Cahn1959a} and the Phase-Field Crystal (PFC) model~\cite{Elder2004,Elder2007,Emmerich2012} are the two most well-known phase-field models, with the former describing phase separation and the latter representing crystallization. The accurate and efficient numerical solution of these equations remains an active area of research, with several established methods tackling the challenges posed by stiff nonlinear partial differential equations (PDEs)~\cite{Biner2017programming}. One of the challenges of phase-field modeling is the presence of multiple length scales, which can lead to stiff problems that are computationally expensive to solve. To overcome numerical difficulties, pseudo-spectral methods have gained attraction because highly accurate and efficient numerical solutions are obtained.

Pseudo-spectral methods~\cite{Boyd2001chebyshev,Shizgal2015spectral} represent a powerful category of numerical methods capable of solving partial differential equations with periodic boundary conditions. These techniques deviate by defining representation basis functions at discrete grid points. This discrete formulation facilitates the efficient evaluation of specific operators, resulting in substantial computational speed-up when employing fast algorithms such as the Fast Fourier Transform (FFT). By discretizing the spatial domain using a Fourier basis and expressing the solution as a sum of Fourier modes, pseudo-spectral methods achieve high accuracy while minimizing computational costs. These methods have been successfully applied to a diverse range of problems in physics and chemistry, including Burgers, Navier-Stokes~\cite{Mortensen2016}, Kardar-Parisi-Zhang~\cite{Giada2002}, Allen-Cahn~\cite{Ayub2019}, Cahn-Hilliard \cite{Yoon2020,Chen1998}, and Phase-Field-Crystal equations~\cite{Tegze2009,Martinez-Agustin2022}. As pseudo-spectral methods have exhibited strong performance in this variety of applications, it is crucial to explore how their combination with advanced time integration schemes can further enhance the solution of stiff nonlinear PDEs.

Recently, three robust numerical methods have gained prominence for solving stiff nonlinear PDEs: Implicit-Explicit (IMEX), Integrating Factor (IF)~\cite{Milewski1999,Fornberg1999,Kassam2005,Krogstad2005}, and Exponential Time Differencing (ETD) methods~\cite{Beylkin1998,Cox2002,Whalen2015}. These exponential integrators have showcased robust performance and accuracy when addressing stiff nonlinear partial differential equations. These methods present an appealing alternative to the widely-used IMEX Euler integrators in phase-field modeling. By capitalizing on the synergy between pseudo-spectral techniques and exponential integrators, we aim to emphasize the significant benefits and improvements in modeling complex systems governed by phase-field dynamics.

Building upon the demonstrated advantages of combining pseudo-spectral techniques with exponential integrators, this paper probes a comprehensive Python implementation of these integrators, showcasing their effectiveness in solving phase-field model equations and emphasizing the potential of Python usage for advancing phase-field modeling. The results obtained from our Python implementation stress the accuracy and computational efficiency of the IMEX and ETD methods, further highlighting the potential of this approach in advancing phase-field modeling research. The codes presented in this work are available online (\url{https://github.com/elvissoares/spectralETD}) under a GPL license.

The paper is organized as follows: In Section~\ref{sec:Theory}, a brief review of the theoretical formalism for phase-field models, pseudo-spectral methods, and time integration schemes is presented. The numerical implementation is detailed in Section~\ref{sec:numerical}. Results of the numerical calculations for phase field models are provided in Section~\ref{sec:results}. Section~\ref{sec:conclusion} concludes this paper.

\section{Formalism}
\label{sec:Theory}

\subsection{Phase-Field Models with Conserved Order Parameters}

Phase-Field models~\cite{Provatas2011} are based on a mathematical concept called a phase field, $\eta(\vb*{r},t)$, an order parameter that is a continuous function that describes the spatial distribution of different phases or microstructures in a material. 

Here we concentrate on the case where $\eta(\vb*{r},t)$ is a conserved order parameter. The conservation equation for $\eta(\vb*{r},t)$ is described by $\pdv*{\eta}{t} + \div{\vb*{j}_\eta} = 0$, with the flux density given by $\vb*{j}_\eta = -M \nabla(\delta F[\eta]/\delta \eta)$, where $F[\eta]$ is the free-energy functional of the phase field. The dynamics of the system can be described by
\begin{align}
    \pdv{\eta}{t} = \div\left[ M(\eta) \nabla \left( \fdv{F[\eta]}{\eta} \right) \right],
    \label{eq:conserved_order_parameter_equation}
\end{align}
where $M(\eta)$ is a positive-defined mobility coefficient that can depend on $\eta$. 

The Cahn-Hilliard (CH) equation is a phase-field model that describes the dynamics of phase separation and coarsening process in two-component mixtures, and it has been applied to various other fields, such as material science, biology, and fluid dynamics. It was introduced in 1958 by John W. Cahn and John E. Hilliard~\cite{Cahn1958} as a mathematical model to understand the thermodynamics and kinetics of phase separation in alloys at a mesoscopic scale, between the atomic and macroscopic scales. The CH free-energy functional is defined as
\begin{align}
    F[\eta(\vb*{r})] = \int_V \dd{\vb*{r}} \left[ \frac{\kappa}{2}\left(\nabla \eta\right)^2 + f(\eta)\right],
    \label{eq:ch_free-energy}
\end{align}
where $\eta(\vb*{r})$ is the concentration field of one phase, $\kappa$ is a positive constant that controls the interfacial energy between the two phases, and the term $f(\eta) = W \eta^2(1-\eta)^2$ defines the double-well free-energy density function, which describes the energetics of the phase separation process. The Eq.~\eqref{eq:conserved_order_parameter_equation} leads to the dynamical equation given by 
\begin{align}
    \pdv{\eta}{t} = M \laplacian \left[ - \kappa \laplacian \eta + 2 W(\eta-3\eta^2+2\eta^3) \right],
    \label{eq:ch_dynamics}
\end{align}
with the $\kappa \laplacian \eta$ term accounting for the interfacial energy. This term tends to smooth out the concentration field and promote the formation of well-defined interfaces between the two phases. Therefore, the CH equation models the time evolution of the order parameter $\eta$ as a result of the competition between the driving force for phase separation, which comes from the free-energy density function $f(\eta)$, and the energy penalty $\kappa \laplacian \eta$ associated with the formation of interfaces between the phases. 

The phase field crystal (PFC) equation describes crystallization growth, where atomic- and microscales are coupled. Like conventional phase field models, the PFC theory incorporates a free energy functional and averages over rapid temporal fluctuations, thereby yielding a time scale of evolution that is on the order of diffusion rather than atomic vibrations. Besides that, the PFC model does not average over atomic distances, leading to the formation of equilibrium patterns that are representative of the crystalline structures.  

Here we focus on the model from Elder and Grant presented in Ref.~\cite{Elder2004}. In this PFC model, the free-energy functional is written as  
\begin{align}
    F[\eta(\vb*{r})] = \int_V \dd{\vb*{r}} \left[ \frac{1}{2} \eta \left(1+\laplacian\right)^2 \eta + \frac{1}{4} \eta^2(2r+\eta^2) \right],
\end{align}
with $r$ being a constant proportional to the temperature deviation from the melting point, and the term $\frac{1}{2} \eta \left(1+\laplacian\right)^2 \eta$ reproduces the first-order peak of the structure factor for the liquid phase.

Using the functional derivative of this free-energy in Eq.~\eqref{eq:conserved_order_parameter_equation}, we get  
\begin{align}
    \pdv{\eta}{t} = M \laplacian \qty[ \nabla^4 \eta + 2\laplacian \eta  + (1+r) \eta + \eta^3],
    \label{eq:pfc_equation}
\end{align}
where the first two terms account for the formation of the interface between the liquid-solid introducing the periodic pattern in the solid phase, respectively. The last two terms in Eq.~\eqref{eq:pfc_equation} constitute a double-well free energy between the liquid and solid phases. 

All the phase field models with conserved order parameters have conservative dynamics that minimize the total free-energy functional, $F[\eta(\vb*{r},t)]$, of the system. In fact, the dynamics described by Eq.~\eqref{eq:conserved_order_parameter_equation} leads to the free-energy evolution given by 
\begin{align}
    \dv{F}{t} = - \int_V \dd{\vb*{r}} M(\eta) \qty[\grad(\fdv{F[\eta]}{\eta})]^2 \geq 0,
    \label{eq:free-energy-evolution}
\end{align}
such that the free energy can just decrease or remain constant during the time evolution of $\eta(\vb*{r},t)$.

\subsection{Pseudo-Spectral Version of Conservative Equations}

The dynamics of any field $\eta(\vb*{r},t)$ described by the Eq.~\eqref{eq:conserved_order_parameter_equation} can be written as the sum of a set of linear terms and another set of non-linear terms in the form 
\begin{align}
    \pdv{\eta}{t} = \sum_{\alpha=0}^A \mathcal{L}_\alpha \nabla^\alpha \eta + \sum_{\beta=0}^B \mathcal{N}_\beta \nabla^\beta f_\beta(\eta),
    \label{eq:general_conservative_nonlinear_equation}
\end{align}
where $\eta$ is the field, $\nabla$ is the spatial coordinate derivative operator, $t$ is time, $\mathcal{L}_\alpha$ and $\mathcal{N}_\beta$ are numerical coefficients of the linear and non-linear terms, respectively, and $f_\beta(\eta)$ are some non-linear functions of the field $\eta(\vb*{r},t)$. Here, $\alpha$ and $\beta$ are nonnegative integer indexes for the $A+1$ linear terms and $B+1$ non-linear terms, respectively. On conservative systems, we expect that $\alpha \geq 2$. 

We can decompose the field $\eta$ into independent normal modes by setting the direct Fourier transform as 
\begin{align}
    \widetilde{\eta}_{\vb*{k}}(t) \equiv \mathcal{F}\{\eta(\vb*{r},t) \}= \int_{L^d} \dd{\vb*{r}} \eta(\vb*{r},t)  e^{-i \vb*{k} \cdot \vb*{r} },
\end{align}
and the inverse Fourier transform as
\begin{align}
    \eta(\vb*{r},t) \equiv \mathcal{F}^{-1}\{\widetilde{\eta}_{\vb*{k}}(t) \}=  \frac{1}{N^d}\sum_{\vb*{k}} \widetilde{\eta}_{\vb*{k}} (t) e^{i \vb*{k} \cdot \vb*{r} },
\end{align}
with $\vb*{k} = \sum_{i=1}^d k_i \vu*{u}_i$ being the wavenumber vector in the Fourier space in $d$ dimensions, and $\mathcal{F}$ and $\mathcal{F}^{-1}$ represent the direct and inverse Fourier transforms, respectively. The complete Fourier space can be defined in each direction by the wavenumbers given by $k_i = \{0,2\pi/L, 4\pi/L,\ldots,2\pi N/L\}$, where $L$ and $N$ are the length of the box and the number of gridpoints in the $i$-direction, respectively.

Any spatial derivative of the field $\eta(\vb*{r},t)$ can be calculated by the Fourier transform as $\nabla^\alpha \eta = \sum_{\vb*{k}} (i \vb*{k})^\alpha \widetilde{\eta}_{\vb*{k}} e^{i \vb*{k} \cdot \vb*{r} }$. Similarly, we have $\nabla^\beta f_\beta(\eta) = \sum_{\vb*{k}} (i \vb*{k})^\beta \mathcal{F}\qty{f_\beta(\eta)}_{\vb*{k}}  e^{i \vb*{k} \cdot \vb*{r} }$. Collecting the Fourier coefficients, we have 
\begin{align}
    \pdv{\widetilde{\eta}_{\vb*{k}}}{t} = \sum_{\alpha=0} \mathcal{L}_\alpha (i \vb*{k})^\alpha \widetilde{\eta}_{\vb*{k}} + \sum_{\beta=0} \mathcal{N}_\beta (i \vb*{k})^\beta \mathcal{F}\qty{f_\beta(\eta)}_{\vb*{k}},
\end{align}
and defining the auxiliary linear and non-linear marching operators, respectively, in the form $\widetilde{\mathcal{L}}_{\vb*{k}} = \sum_{\alpha=0} \mathcal{L}_\alpha (i \vb*{k})^\alpha,$ and $\widetilde{\mathcal{N}}_{\vb*{k}} (\eta(t)) = \sum_{\beta=0} \mathcal{N}_\beta (i \vb*{k})^\beta \mathcal{F}\qty{f_\beta(\eta)}_{\vb*{k}}$, we can simplify the FT dynamical equation as 
\begin{align}
    \pdv{\widetilde{\eta}_{\vb*{k}}}{t} = \widetilde{\mathcal{L}}_{\vb*{k}} \widetilde{\eta}_{\vb*{k}} + \widetilde{\mathcal{N}}_{\vb*{k}} (\eta(t)).
    \label{eq:general_model_ft}
\end{align}
Note that the linear marching operator is time-independent while the non-linear marching operator is time-dependent because the non-linear functions $f_\beta$ are dependent on $\eta(t)$. 

If there aren't any non-linear terms, \emph{i.e.}, $f_\beta= 0$, the equation has an analytical solution in the form $\widetilde{\eta}_{\vb*{k}}(t) = \widetilde{\eta}_{\vb*{k}}(0) e^{t \widetilde{\mathcal{L}}_{\vb*{k}}}$. Therefore, the product $\widetilde{\eta}_{\vb*{k}}(t)e^{-t \widetilde{\mathcal{L}}_{\vb*{k}}}$ is a time-invariant of the pure linear model. In this manner, the exponential term $e^{-t \widetilde{\mathcal{L}}_{\vb*{k}}}$ can be used as an integrand factor to solve the more general problem, Eq.~\eqref{eq:general_model_ft}. Multiplying the Eq.~\eqref{eq:general_model_ft} by the integrand factor $e^{-t \widetilde{\mathcal{L}}_{\vb*{k}}}$, we can identify 
\begin{align}
    \pdv{}{t}\qty(\widetilde{\eta}_{\vb*{k}}(t)e^{-t \widetilde{\mathcal{L}}_{\vb*{k}}}) = \widetilde{\mathcal{N}}_{\vb*{k}} (\eta(t)) e^{-t \widetilde{\mathcal{L}}_{\vb*{k}}},
    \label{eq:integrating_factor}
\end{align}
without any loss of generality. Integrating Eq.~\eqref{eq:integrating_factor} from $t_0$ to $t$, we get 
\begin{align}
    \widetilde{\eta}_{\vb*{k}}(t) = \qty[\widetilde{\eta}_{\vb*{k}}(t_0) + \int_0^{t-t_0}\widetilde{\mathcal{N}}_{\vb*{k}} (\eta(\tau+t_0)) e^{-\tau \widetilde{\mathcal{L}}_{\vb*{k}}} \dd{\tau}]e^{(t-t_0) \widetilde{\mathcal{L}}_{\vb*{k}}}.
    \label{eq:general_formula}
\end{align}
which is an analytical relation without any approximations.

\subsection{Time Integration Methods}
\label{sec:time_integration}

\subsubsection{Implicit-Explicit Euler Method: IMEX}

The IMEX (IM = implicit, EX = explicit) method \cite{Chen1998} can be used to approximate the time integration of Eq.~\eqref{eq:general_model_ft}. The method can be obtained by discretizing Eq.~\eqref{eq:general_model_ft} from $t_n$ until $t_{n+1} = t_n + h$, using an implicit step for the linear term and an explicit step for the nonlinear term, as follows
\begin{align}
    \frac{\widetilde{\eta}_{\vb*{k}}(t_{n+1}) - \widetilde{\eta}_{\vb*{k}}(t_n)}{h} = \widetilde{\mathcal{L}}_{\vb*{k}} \widetilde{\eta}_{\vb*{k}}(t_{n+1})+ \widetilde{\mathcal{N}}_{\vb*{k}} (\eta(t_n)),
    \label{eq:deduced_SIE}
\end{align}
where $t_{n+1} = t_n+ h$. The Eq.~\eqref{eq:deduced_SIE} can be solved for $\widetilde{\eta}_{\vb*{k}}(t_{n+1})$ such that 
\begin{align}
    \widetilde{\eta}_{\vb*{k}}(t_{n+1}) = \frac{\widetilde{\eta}_{\vb*{k}}(t_n)+ h \widetilde{\mathcal{N}}_{\vb*{k}} (\eta(t_n))}{1-h \widetilde{\mathcal{L}}_{\vb*{k}}}.
\end{align}

The implicit scheme is stable, so larger time steps may be used. But there is a necessary condition of $h |\widetilde{\mathcal{L}}_{\vb*{k}}| < 1.0$ to make the convergence. 

\subsubsection{Integrating Factor Method: IF}

The IF method~\cite{Milewski1999,Fornberg1999,Kassam2005,Krogstad2005} can be obtained by discretizing the Eq.~\eqref{eq:integrating_factor} in time from $t_n$ until $t_{n+1} = t_n + h$, using an explicit step for the nonlinear source term in the following form 
\begin{align}
    \frac{\widetilde{\eta}_{\vb*{k}}(t_{n+1})e^{-t_{n+1} \widetilde{\mathcal{L}}_{\vb*{k}}} - \widetilde{\eta}_{\vb*{k}}(t_n)e^{-t_n \widetilde{\mathcal{L}}_{\vb*{k}}}}{h}  = \widetilde{\mathcal{N}}_{\vb*{k}} (\eta(t_n)) e^{-t_n \widetilde{\mathcal{L}}_{\vb*{k}}},
\end{align}
such that 
\begin{align}
    \widetilde{\eta}_{\vb*{k}}(t_{n+1})  = \widetilde{\eta}_{\vb*{k}}(t_n)e^{h \widetilde{\mathcal{L}}_{\vb*{k}}} + h \widetilde{\mathcal{N}}_{\vb*{k}} (\eta(t_n)) e^{h \widetilde{\mathcal{L}}_{\vb*{k}}}.
    \label{eq:IF}
\end{align}

If we multiply Eq.~\eqref{eq:IF} by $e^{-h\widetilde{\mathcal{L}}_{\vb*{k}}}$, and take the limit $h|\widetilde{\mathcal{L}}_{\vb*{k}}| \ll 1$ such that $e^{-h\widetilde{\mathcal{L}}_{\vb*{k}}} \approx 1 - h\widetilde{\mathcal{L}}_{\vb*{k}}$, we obtain the IMEX method. 

\subsubsection{Exponential Time Differencing Method: ETD}

The ETD method~\cite{Pope1963,Cox2002} consists of discretizing Eq.~\eqref{eq:general_formula} in time from $t_n$ until $t_{n+1} = t_n + h$. Using a simple constant approximation of the non-linear operator in the integral $\int_0^{h}\widetilde{\mathcal{N}}_{\vb*{k}} (\eta(\tau+t_n)) e^{-\tau \widetilde{\mathcal{L}}_{\vb*{k}}} \dd{\tau} \approx  \widetilde{\mathcal{N}}_{\vb*{k}}(\eta(t_n))[1-e^{-h \widetilde{\mathcal{L}}_{\vb*{k}}}]/\widetilde{\mathcal{L}}_{\vb*{k}}$, such that 
\begin{align}
    \widetilde{\eta}_{\vb*{k}}(t_{n+1}) = \widetilde{\eta}_{\vb*{k}}(t_n)e^{h  \widetilde{\mathcal{L}}_{\vb*{k}}} + h \widetilde{\mathcal{N}}_{\vb*{k}}(\eta(t_n))\qty[\frac{e^{h  \widetilde{\mathcal{L}}_{\vb*{k}}}-1}{h \widetilde{\mathcal{L}}_{\vb*{k}}}].
    \label{eq:ETD}
\end{align}

Note that, multiplying Eq.~\eqref{eq:ETD} by $e^{-h\widetilde{\mathcal{L}}_{\vb*{k}}}$ and taking $|\widetilde{\mathcal{L}}_{\vb*{k}}| \ll 1$, we obtain the IMEX method. 

The function $(e^x - 1)/x$, despite being continuous, presents undesirable behavior near $x = 0$ in a computational context due to cancellation errors. This can be bypassed by using the appropriate limit of the function at $x=0$.

\section{Numerical Implementation in Python}
\label{sec:numerical}

We have implemented all the time integration schemes combined with the pseudo-spectral method, as described in the previous section, in high-level Python code.

For the 1D geometry, we comment on the cases of the advection-diffusion equation and the Burgers' equation in Appendix~\ref{app:other}.

In the examples considered here, we consider a 2D geometry as our physical space that consists of a structured uniform square $[0, 16\pi]^2$, where each direction is divided into $N=2^M$ uniform intervals, where $M$ is a positive integer. The computational mesh is represented in Python as
\begin{lstlisting}[language=Python]
import numpy as np
# Size of the system
N = 2**8 # 2**8 = 256
L = 16*np.pi
x = np.linspace(0,L,N)
dx = x[1]-x[0]
\end{lstlisting}

The system was evolved with a stepsize of $h = 0.01$, and simulation was carried out up to 150 000 time steps. At each 100th frame, we output and store the result in an array. 
\begin{lstlisting}[language=Python]
# The time step definition
h = 0.01
T = 1500
Nsteps = int(T/h)
dframes = 1.0 # time step to output
Nframes = int(T/dframes) #frames to the output
nframes = Nsteps//Nframes
# The array of outputs
n = np.empty((Nframes,N,N), dtype=np.float32)
\end{lstlisting}

Our pseudo-spectral Python solver extensively uses the NumPy~\cite{harris2020array} package. The FFT module comes from Scipy~\cite{2020SciPy} package. The fourier space can be defined by 
\begin{lstlisting}[language=Python]
# The Fourier variables
from scipy.fft import fft2, ifft2
n_k = np.empty((N,N), dtype=np.complex64)
kx = np.fft.fftfreq(N, d=dx)*2*np.pi
k = np.array(np.meshgrid(kx , kx ,indexing ='ij'), dtype=np.float32)
k2 = np.sum(k*k,axis=0, dtype=np.float32)
\end{lstlisting}

To avoid any aliasing problem with non-linear terms due to the discrete nature of numerical simulations, we have used the \emph{2/3 rule} dealising technique~\cite{Orszag1971}. The aliasing happens when the signal or function is not adequately sampled, violating the Nyquist-Shannon sampling theorem. Here, we introduce the dealising matrix as $\xi^\text{dealias}_{\vb*{k}} =  \Pi_i^d \Theta\qty(k_i^\text{cut}-k_i)$, where $k_i^\text{cut}=(2/3)k_i^\text{max}$ is the Nyquist cutoff frequency in the $i$th direction and $\Theta(x)$ is the Heaviside step function. 

\begin{lstlisting}[language=Python]
kmax_dealias = kx.max()*2.0/3.0 # The Nyquist mode
# Dealising matrix
dealias = np.array((np.abs(k[0]) < kmax_dealias )*(np.abs(k[1]) < kmax_dealias ),dtype =bool)
\end{lstlisting}

The linear and nonlinear marching operators can be defined by 
\begin{lstlisting}[language=Python]
# The linear terms of PDE
Loperator_k = (...) # some function of k only
# The non-linear terms of PDE 
def Noperator_func(n):
    return (...) # some function of n and k
\end{lstlisting}

After that, the time integration schemes can be chosen using the following code
\begin{lstlisting}[language=Python]
# Defining the time marching operators arrays
# can be calculated once
if method == 'IMEX':
    Tlinear_k = 1.0/(1.0-h*Loperator_k) 
    Tnon_k = dealias*h/(1.0-h*Loperator_k) 
elif method == 'IF':
    Tlinear_k = np.exp(h*Loperator_k) 
    Tnon_k = dealias*h*Tlinear_k
elif method == 'ETD':
    Tlinear_k = np.exp(h*Loperator_k) 
    def myexp(x):
        if x == 1: return 1.0
        else: return (x-1.0)/np.log(x)
    vmyexp = np.vectorize(myexp) # vectorize myexp
    Tnon_k = dealias*h*vmyexp(Tlinear_k)
else: print('ERROR: Undefined Integrator')
\end{lstlisting}

Note that the operator can be calculated once and before the time integration loop. 

Finally, the loop of time integration is 
\begin{lstlisting}[language=Python]
n[0] = (...) # some initial condition
Noperator_k = n_k.copy() # auxiliary array
nn = n[0].copy() # auxiliary array
n_k[:] = fft2(n[0]) # FT initial condition
# time evolution loop
for i in range(1,Nsteps):
    # calculate the nonlinear operator
    Noperator_k[:] = Noperator_func(nn)
    # updating in time
    n_k[:] = n_k*Tlinear_k + Noperator_k*Tnon_k 
    # IFT to next step
    nn[:] = ifft2(n_k).real 
    # test to output
    if (i % nframes) == 0: n[i//nframes] = nn
\end{lstlisting}

\section{Results and Discussion} 
\label{sec:results}

\subsection{2D Cahn-Hilliard Equation}

Here we consider the CH equation with the parameters set as $W=1$, $\kappa=0.1$, and $M=1$. The initial profile is a homogeneous field with value $\eta_0$ and a Gaussian perturbation of intensity 0.02 at each gridpoint. 
\begin{lstlisting}[language=Python]
# Cahn-Hilliard model constants
W = 1.0
M = 1.0 # mobility
kappa = 0.1 #gradient coeficient
# Initial condition
rng = np.random.default_rng(12345) # the seed
noise = 0.02
n0 = 0.5
n[0] = n0 +noise*rng.standard_normal(n[0].shape)
\end{lstlisting}

From the pseudo-spectral method, the CH equation of motion, Eq.~\eqref{eq:ch_dynamics}, can be re-written as Eq.~\eqref{eq:general_model_ft} with the time-marching operators given by 
\begin{align}
    \widetilde{\mathcal{L}}_{\vb*{k}} &= - M (\kappa k^4 +2 W k^2), \\
    \widetilde{\mathcal{N}}_{\vb*{k}}(\eta) &= - 2M  W k^2 \mathcal{F}\qty{-3\eta^2+2\eta^3}_{\vb*{k}},
\end{align}
with the following Python code: 
\begin{lstlisting}[language=Python]
# The linear terms of PDE
Loperator_k = -M*(kappa*k2**2+2*W*k2)
# The non-linear terms of PDE 
def Noperator_func(n):
    return -2*M*W*k2*fft2(-3*n**2+2*n**3)
\end{lstlisting}

The whole Python code to solve the CH equation with exponential integrators and pseudo-spectral methods is presented in Appendix~\ref{app:python}. 

Figure~\ref{fig:ch} shows the result of the CH equation using the pseudo-spectral method with the ETD integration scheme for three different initial condition $\eta_0=0.4$, 0.5 and 0.6 (see Movies S1, S2, and S3 for more details). The first row presents the dynamics of the nucleation and growth of bubbles of the red phase inside the continuous  blue phase. The second row represents the dynamics of spinodal decomposition of the two phases. And the third row presents the nucleation dynamics of the blue phase surrounded by the continuous red phase. 

\begin{figure*}[htbp]
    \centering
    \includegraphics{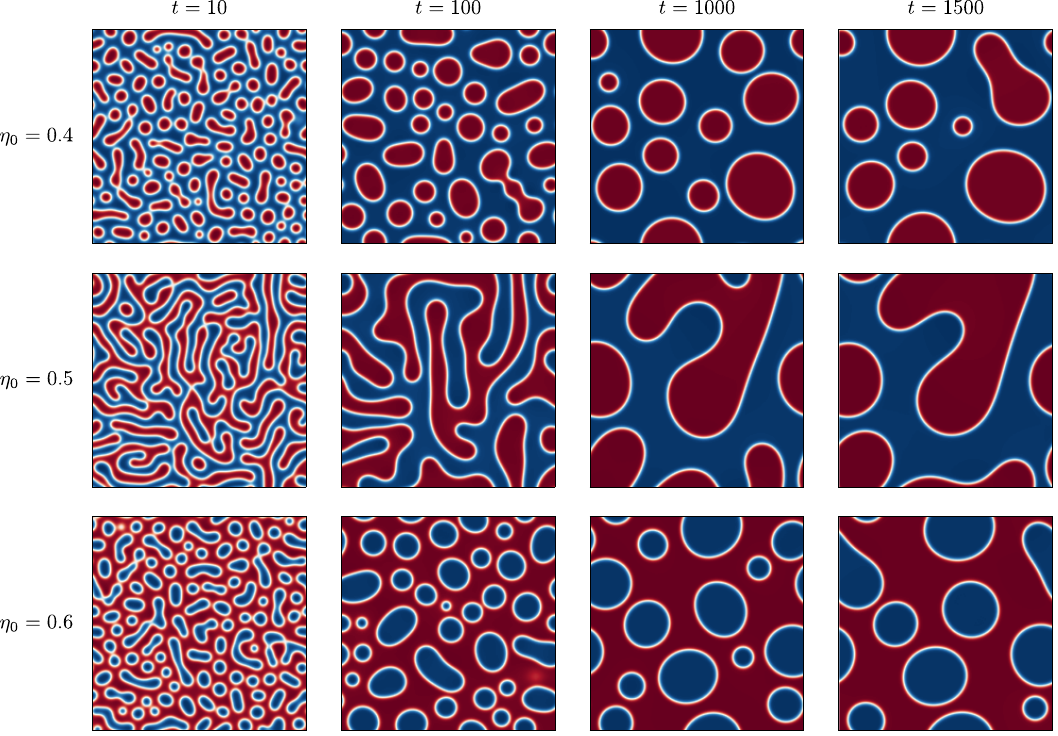}
    \caption{The time evolution of the field profile $\eta(\vb*{r},t)$ calculated from the CH equation. The different rows represent different initial conditions: $\eta_0 = 0.4$ (see Movie S1), $\eta_0 = 0.5$ (see Movie S2), and $\eta_0 = 0.6$ (see Movie S3). The columns represents the snapshots in different instant of time $t$. The system was evolved with the ETD scheme and a stepsize of $h = 0.01$. The initial profile is a homogeneous field with $\eta_0$ and a Gaussian perturbation of intensity 0.02, \emph{i.e.}, $\eta(t=0) = \eta_0 + 0.02 \mathcal{N}(0,1)$. }
    \label{fig:ch}
\end{figure*} 

\begin{figure}[htbp]
    \centering
    \includegraphics{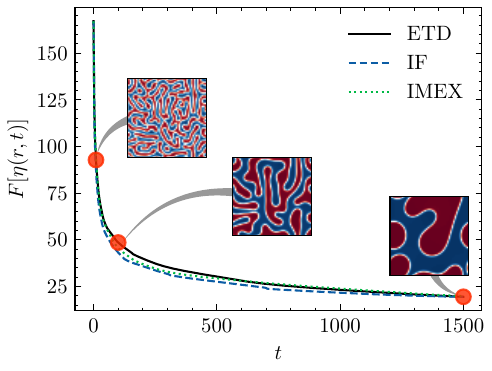}
    \caption{Free energy as a function of time to the CH equation with the initial condition as a homogeneous field with $\eta_0=0.5$ and a Gaussian perturbation of intensity 0.02 with integrating using different methods with a timestep of $h=0.01$. The solid line represents the solution with the ETD method, the dashed line with the IF method, and the dotted line with the IMEX method.}
    \label{fig:freeenergy-ch}
\end{figure}

Figure~\ref{fig:freeenergy-ch} illustrates the temporal evolution of the CH free energy, as defined by Eq.~\eqref{eq:ch_free-energy}, during spinodal decomposition ($\eta_0=0.5$) for the three distinct time integration schemes. It can be observed intense decay in the free-energy during the initial times up to $t=100$. Following that, the system relaxes at a slower rate for approximately 1400 units of time. The various integration schemes display minor variations in the free energy decay curve over the course of the process, although they converge to a similar energy value towards the end of the numerical calculation.

\begin{table}[htbp]
    \centering
    \caption{The $L_2$-error of the numerical solution of the CH equation for the spinodal decomposition condition with $\eta_0=0.5$ at $t=100$.}
    \begin{ruledtabular}
    \begin{tabular}{cccc}
    $h$ & IMEX & IF & ETD \\ \hline
    0.0002 & $6.4 \times 10^{-4}$ & $9.4 \times 10^{-4}$ & $6.4 \times 10^{-4}$ \\
    0.001 & $7.4 \times 10^{-4}$ & $1.3 \times 10^{-3}$ & $6.1 \times 10^{-4}$ \\
    0.002 & $7.5 \times 10^{-4}$ & $1.7 \times 10^{-3}$ & $7.4 \times 10^{-4}$ \\
    0.01 & $1.4 \times 10^{-3}$ & $2.2 \times 10^{-3}$ & $1.3 \times 10^{-3}$ \\
    0.02 & $1.8 \times 10^{-3}$ & $2.4 \times 10^{-3}$ & $1.4 \times 10^{-3}$ \\
    0.1 & $2.0 \times 10^{-3}$ & $2.4 \times 10^{-3}$ & $2.0 \times 10^{-3}$ \\
    0.2 & $2.2 \times 10^{-3}$ & $2.2 \times 10^{-3}$ & $2.2 \times 10^{-3}$ \\
    1.0 & $2.3 \times 10^{-3}$ & $1.9 \times 10^{-3}$ & $2.4 \times 10^{-3}$
    \end{tabular}
    \end{ruledtabular}
    \label{tab:error_CH}
\end{table} 

The $L_2$-error is defined as the absolute error given by $ \norm{\eta - \eta_r}_2/N^2$ where $\eta$ is the numerical solution with determined stepsize $h$ and determined scheme, and $\eta_r$ is the numerical solution with the reference stepsize of $10^{-4}$ and ETD scheme. The normalization $N^2$ is the total number of gridpoints. Table~\ref{tab:error_CH} presents the $L_2$-error of the three different time integration schemes. The ETD method comes up as the most accurate when compared to the IF method. However, the simple IMEX method remains competitive even when taking larger time steps into account.

\subsection{2D Phase Field Crystal Equation}

We consider the PFC equation with the parameters set as $r=-0.25$ and $M=1.0$. The initial profile is a homogeneous field $\eta_0$ with a Gaussian perturbation of intensity $0.02|\eta_0|$ at each gridpoint. 
\begin{lstlisting}[language=Python]
# PFC model constants
r = -0.25
M = 1.0 # mobility
# Initial condition
rng = np.random.default_rng(12345)
n0 = -0.285
noise = 0.02*np.abs(n0)
n[0] = n0 +noise*rng.standard_normal(n[0].shape)
\end{lstlisting}

From the pseudo-spectral method, the PFC equation can be re-written as Eq.~\eqref{eq:general_model_ft} with the time-marching operators given by 
\begin{align}
    \widetilde{\mathcal{L}}_{\vb*{k}} &= - M k^2(k^4 - 2k^2 + 1 + r),\\
    \widetilde{\mathcal{N}}_{\vb*{k}}(\eta) &= - M k^2 \mathcal{F}\qty{\eta^3}_{\vb*{k}},
\end{align}
with the following Python code:
\begin{lstlisting}[language=Python]
# The linear terms of PDE
Loperator_k = -M*k2*(k2**2-2*k2+1+r)
# The non-linear terms of PDE (with dealising)
def Noperator_func(n):
    return -(k2*M*fft2(n**3))
\end{lstlisting}

Again, the whole Python code to solve the PFC equation with exponential integrators and pseudo-spectral methods is presented in Appendix~\ref{app:python}. 

Figure~\ref{fig:pfc} presents the time evolution of the field profile calculated from the PFC equation using the pseudo-spectral method with the ETD integration scheme. The first line presents the snapshots of the temporal evolution of a system with an initial condition of $\eta_0 = -0.085$, which leads to the formation of a lamellar phase. The second line presents the snapshots of the system's temporal evolution with $\eta_0 = -0.285$, which results in an ordered crystalline phase. (see Movies S4 and S5 for more details).

\begin{figure*}[htbp]
    \centering
    \includegraphics{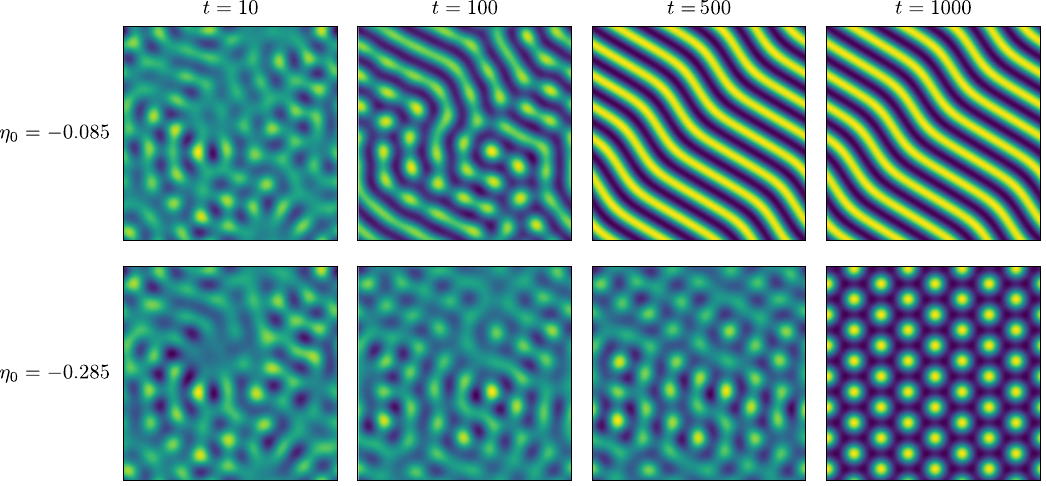}
    \caption{The time evolution of the field profile $\eta(\vb*{r},t)$ calculated from the PFC equation. The different rows represent the two different initial conditions characterized by $\eta_0 = -0.085$ (see Movie S4) and $\eta_0 = -0.085$ (see Movie S5). The columns represent the snapshots of different time instant $t$. The system was evolved with an ETD scheme and stepsize of $h = 0.01$. The initial profile is a homogeneous field with $\eta_0$ and a Gaussian perturbation of intensity $0.02|\eta_0|$, \emph{i.e.}, $\eta(t=0) = \eta_0 + 0.02|\eta_0| \mathcal{N}(0,1)$. }
    \label{fig:pfc}
\end{figure*}  

\begin{figure}[htbp]
    \centering
    \includegraphics{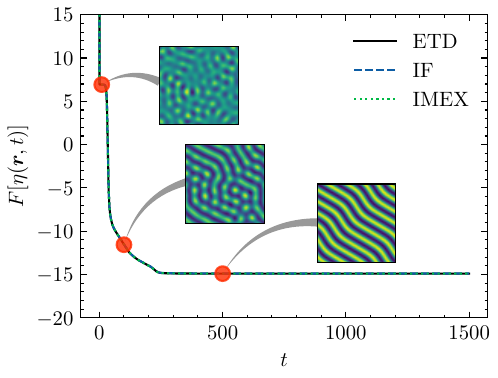}
    \caption{Free energy as a function of time from the PFC equation with the initial condition as a homogeneous field with $\eta_0=-0.085$ and a Gaussian perturbation of intensity $0.02|\eta_0|$.  }
    \label{fig:freeenergy_lamellar}
\end{figure}

Figure~\ref{fig:freeenergy_lamellar} represents the temporal evolution of the free energy for the system generating a lamellar phase derived from the PFC equation. A significant topological transition occurs around $t=10$, leading to the initial formation of a lamellar phase. Between the time interval from $t=100$ to $t=300$, the lamellar phase forms throughout the entire system. The free energy remains constant after this period until the final time step.
  
\begin{figure}[htbp]
    \centering
    \includegraphics{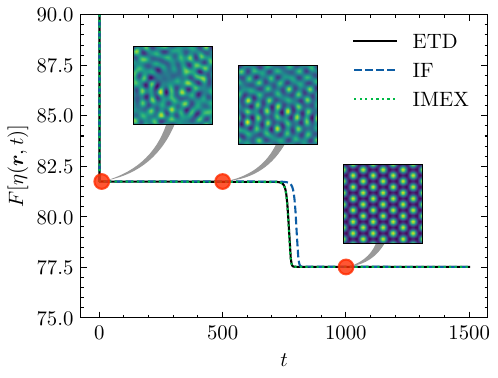}
    \caption{Free energy as a function of time from the PFC equation with the initial condition as a homogeneous field with $\eta_0=-0.285$ and a Gaussian perturbation of intensity $0.02|\eta_0|$. }
    \label{fig:freeenergy_crystal}
\end{figure}

Figure~\ref{fig:freeenergy_crystal} presents the free-energy as a function of time during the crystal formation on the PFC dynamics. There is a sharp decay of the free-energy until $t=10$ due to the decay of the high energy modes. Between $t=10$ and 700, the system experiences a slight variation in free energy, with minimal movements in the phase field configuration. However, near $t=750$, there is a rapid decrease in free energy during the topological transition to form the ordered crystalline lattice. Indeed, $t=750$ is the time required for the system to form the crystal and is associated with the system's crystalline nucleation rate. It is interesting to note that this jump in free energy only occurs during the formation of crystalline phases and does not occur during the formation of the lamellar phase, according to Fig.~\ref{fig:freeenergy_lamellar}.

\begin{table}[htbp]
    \centering
    \caption{The $L_2$-error of the numerical solution of PFC equation for the crystallization condition with $\eta_0=-0.285$ at $t=750$.}
    \begin{ruledtabular}
    \begin{tabular}{cccc}
    $h$ & IMEX & IF & ETD \\ \hline
    0.0002 & $1.5\times 10^{-6}$ & $1.5\times 10^{-6}$ & $1.8\times 10^{-6}$ \\
    0.001 & $5.5 \times 10^{-6}$ & $5.4 \times 10^{-6}$ & $7.0 \times 10^{-6}$ \\
    0.002 & $5.5 \times 10^{-6}$ & $2.3 \times 10^{-6}$ & $5.5 \times 10^{-6}$ \\
    0.01 & $5.1 \times 10^{-6}$ & $1.2 \times 10^{-5}$ & $4.8 \times 10^{-6}$ \\
    0.02 & $5.3 \times 10^{-6}$ & $3.1 \times 10^{-5}$ & $4.8 \times 10^{-6}$ \\
    0.1 & $9.5 \times 10^{-6}$ & $1.8 \times 10^{-4}$ & $6.9 \times 10^{-6}$ \\
    0.2 & $1.4 \times 10^{-5}$ & $2.3 \times 10^{-4}$ & $9.4 \times 10^{-6}$ \\
    1.0 & $3.0 \times 10^{-5}$ & $2.3 \times 10^{-4}$ & $2.3 \times 10^{-5}$
    \end{tabular}
    \end{ruledtabular}
    \label{tab:error_PFC}
\end{table} 

Table~\ref{tab:error_PFC} presents the $L_2$-error of the three different time integration schemes in the case of the PFC equation for crystalization. The ETD method comes up as the most accurate when compared to the other two methods. The ETD method remains very precise even when taking larger time steps into account. It is worth noting that the IF method leads to a different nucleation rate when compared to the other two methods. Therefore, avoiding the IF method when dealing with PFC equations is interesting.

\section{Conclusions}
\label{sec:conclusion}

A comprehensive investigation has been conducted, focusing on three distinct time integration schemes (IMEX, IF, and ETD) in combination with pseudo-spectral methods, which are designed to lead with stiff or highly oscillatory nonlinear PDEs, thereby highlighting the potential utility of these methodologies in addressing complex scientific problems. 

The ETD method demonstrates its superiority in accuracy and enhanced stability compared to the IF method, solidifying its position as a preferred choice for handling complex numerical problems. Moreover, avoiding the IF method when addressing PFC equations is important, as it results in a distinct nucleation rate compared to the other two methods. However, it is worth noting that the simple IMEX method maintains its competitiveness. This suggests that, depending on the specific requirements and constraints of a given problem, the IMEX method may still be a viable alternative for first applications. After all, the ETD method exhibits enhanced stability compared to other methods, further contributing to its appeal as a robust and reliable choice for tackling various phase-field modeling scenarios. In summary, weighing the trade-offs between accuracy, stability, and computational efficiency is essential when selecting the most appropriate method for a particular phase-field modeling scenario.

The combination of pseudo-spectral techniques and exponential integrators showcases significant benefits for modeling complex systems governed by phase-field dynamics, including solidification processes and pattern formation. Our comprehensive Python implementation shows the effectiveness of this combined approach in solving phase-field model equations, highlighting the accuracy and computational advantages of the ETD method compared to other numerical approaches.

In conclusion, this work underscores the potential of leveraging exponential integrators and pseudo-spectral techniques to advance phase-field modeling research and provide accurate and efficient solutions for complex scientific problems. We aim to expand this methodology by improving ETD via Runge-Kutta schemes, as already discussed in Ref.~\cite{Cox2002} for other equations. We are also interested in examining its applicability in the context of stochastic partial differential equations.

\section*{Acknowledgments}

The authors thank Petrobras and Shell, which provided financial support through the Research, Development, and Innovation Investment Clause, in collaboration with the Brazilian National Agency of Petroleum, Natural Gas, and Biofuels (ANP, Brazil). Additionally, this research was partially funded by CNPq, CAPES, and FAPERJ.

\section*{The data and code availability}
The data and code that support the findings of this study are available in the author's Github repository: \url{https://github.com/elvissoares/spectralETD}.

\section*{Supplementary Material}
Supplementary material for this article is available at: \textcolor{red}{link}
\begin{itemize}
    \item \emph{Movie S1}. Time evolution of 2D CH equation ($W=1$, $\kappa=0.1$ and $M=1$) with $\eta_0 = 0.4$ (red phase nucleation-growth-ripening) using pseudo-spectral method with ETD scheme.
    \item \emph{Movie S2}. Time evolution of 2D CH equation ($W=1$, $\kappa=0.1$ and $M=1$) with $\eta_0 = 0.5$ (spinodal decomposition) using pseudo-spectral method with ETD scheme.
    \item \emph{Movie S3}. Time evolution of 2D CH equation ($W=1$, $\kappa=0.1$ and $M=1$) with $\eta_0 = 0.6$ (blue phase nucleation-growth-ripening) using pseudo-spectral method with ETD scheme.
    \item \emph{Movie S4}. Time evolution of 2D PFC equation ($r=-0.25$ and $M=1$) with $\eta_0 = -0.085$ (lamellar phase) using pseudo-spectral method with ETD scheme.
    \item \emph{Movie S5}. Time evolution of 2D PFC equation ($r=-0.25$ and $M=1$) with $\eta_0 = -0.285$ (crystal phase) using pseudo-spectral method with ETD scheme.
    \item \emph{Movie S6}. Time evolution of 1D advection-diffusion equation ($u = 5$ and $D=0.01$) using the analytical solution in Fourier space.
    \item \emph{Movie S7}. Time evolution of 1D Burgers' equation ($\nu=0.001$) using pseudo-spectral method with ETD scheme.
\end{itemize}

\appendix

\section{Other Examples}
\label{app:other}

\subsection{1D Advection-Diffusion Equation}

The advection-diffusion equation is a fundamental partial differential equation that describes the transport of a scalar quantity (e.g., temperature, concentration) in a fluid flow. It combines the effects of advection, representing the transport of the scalar quantity by the fluid flow, and diffusion, which represents the random movement of the scalar quantity due to molecular motion. The one-dimensional advection-diffusion equation is written as
\begin{align}
    \frac{\partial \eta}{\partial t} = - u\frac{\partial \eta}{\partial x} + D \frac{\partial^2 \eta}{\partial x^2},
\end{align}
with $u$ being the velocity constant, $D$ being the diffusion coefficient, and $\eta$ represents the concentration field. The time derivative of the Fourier Transform $\widetilde{\eta}_{\vb*{k}}$ is 
\begin{align}
    \pdv{}{t}\widetilde{\eta}_{\vb*{k}}(t) = -(i u k + D k^2)\widetilde{\eta}_{\vb*{k}}(t).
\end{align}

This case corresponds to the pseudo-spectral equation Eq.~\eqref{eq:general_model_ft} with the following operators given by 
\begin{align}
    \widetilde{\mathcal{L}}_{\vb*{k}} = -(i u k + D k^2), \quad \text{and} \quad \widetilde{\mathcal{N}}_{\vb*{k}} = 0,
\end{align}
such that the time evolution has an analytical solution in the form
\begin{align}
    \widetilde{\eta}_{k}(t) = \widetilde{\eta}_{k}(0) e^{-(i u k + D k^2) t}.
\end{align}

As an example, we set $u = 5$ and $D=0.01$. The size of the system is $L=2\pi$ with the number of gridpoints $N=2^{12} = 4096$. The initial profile is a top-hat function with a length of $l=0.2$ and intensity of $\eta_0 = 1.0$ starting at $x_0 = -\pi$.

\begin{figure}[htbp]
    \centering
    \includegraphics{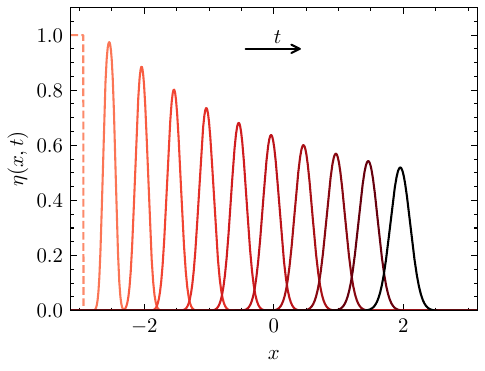}
    \caption{The time evolution of the diffusion-advection equation profile with $D=0.01$ and $u = 5$. The different colored solid lines represent $t=0,0.1,0.2,\ldots,1.0$  with stepsize of $h = 0.001$. The initial profile is a top-hat function represented by the dashed line.}
    \label{fig:diffusion}
\end{figure}

Figure~\ref{fig:diffusion} illustrates the temporal evolution of the system governed by the Advection-Diffusion equation. The initial condition is represented by the dashed line profile. The advection term transports the system from left to right, while the diffusion term smooths the initial distribution, causing it to approach a Gaussian distribution more closely.

\subsection{1D Burgers' Equation}

Burgers' equation describes the dynamics of a fluid with nonlinear advection and diffusion being commonly used as a model for a range of physical phenomena, including turbulence, shock waves, and traffic flow. Here, we apply our approach to solve the same problem as an example of the algorithm. The Burgers' equation is defined as 
\begin{align}
    \frac{\partial \eta}{\partial t} =- \eta\frac{\partial \eta}{\partial x} + \nu \frac{\partial^2 \eta}{\partial x^2}
    \label{eq:burgers}
\end{align}
with $\nu$ being the viscosity coefficient and $\eta$ being the velocity field. The nonlinear term can be rewritten as $\eta\pdv*{\eta}{x}= \tfrac{1}{2}\pdv*{\eta^2}{x}$ such that the pseudo-spectral formulation is defined by Eq.~\eqref{eq:general_model_ft} with the operators written as
\begin{align}
    \widetilde{\mathcal{L}}_{\vb*{k}} = - \nu k^2 \quad \text{and} \quad
    \widetilde{\mathcal{N}}_{\vb*{k}}(\eta) = - \frac{1}{2} i k \mathcal{F}\qty{\eta^2}_{\vb*{k}}.
\end{align}

As an example, we set the viscosity coefficient as $\nu = 0.001$. The size of the system is $L=2\pi$ with the number of gridpoints $N=2^{12}=4096$. The initial profile is a Gaussian given by $\eta_0(x) = \exp(-10(x/2)^2)$.

\begin{figure}[htbp]
    \centering
    \includegraphics{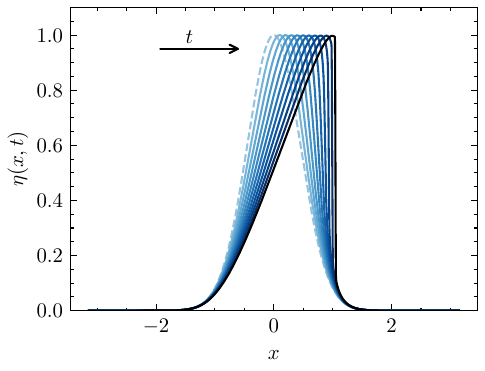}
    \caption{The time evolution of the Burgers' equation profile with $\nu=0.001$. The different colored solid lines represent $t=0,0.1,0.2,\ldots,1.0$ with stepsize of $h = 0.001$. The initial profile is a Gaussian function represented by the dashed line.}
    \label{fig:burgers}
\end{figure}

Figure~\ref{fig:burgers} illustrates the time evolution of the $\eta$ field profile of the Burgers' equation. We can see the appearance of a shock front in the late stages. As these shock waves propagate through the medium, they gradually dissipate due to the  viscosity effects.

\section{Python Codes}
\label{app:python}

\subsection{2D Cahn-Hilliard Model}

\begin{lstlisting}[language=Python]
import numpy as np
from scipy.fft import fft2, ifft2
# Cahn-Hilliard model constants
W = 1.0
M = 1.0 # mobility
kappa = 0.1 #gradient coeficient
# Size of the system
N = 2**8 # 2**8 = 256
L = 16*np.pi
x = np.linspace(0,L,N)
dx = x[1]-x[0]
# The time step definition
h = 0.01
T = 1500
Nsteps = int(T/h)
dframes = 1.0 # time step to output
Nframes = int(T/dframes) #frames to the output
nframes = Nsteps//Nframes
# The array of outputs
n = np.empty((Nframes,N,N), dtype=np.float32)
# The Fourier variables
n_k = np.empty((N,N), dtype=np.complex64)
kx = np.fft.fftfreq(N, d=dx)*2*np.pi
k = np.array(np.meshgrid(kx , kx ,indexing ='ij'), dtype=np.float32)
k2 = np.sum(k*k,axis=0, dtype=np.float32)
kmax_dealias = kx.max()*2.0/3.0 # The Nyquist mode
# Dealising matrix
dealias = np.array((np.abs(k[0]) < kmax_dealias )*(np.abs(k[1]) < kmax_dealias ),dtype =bool)
# The linear terms of PDE
Loperator_k = -M*(kappa*k2**2+2*W*k2)
# The non-linear terms of PDE 
def Noperator_func(n):
    return -2*M*W*k2*fft2(-3*n**2+2*n**3)
# Defining the time marching operators arrays
# can be calculated once
if method == 'IMEX':
    Tlinear_k = 1.0/(1.0-h*Loperator_k) 
    Tnon_k = dealias*h/(1.0-h*Loperator_k) 
elif method == 'IF':
    Tlinear_k = np.exp(h*Loperator_k) 
    Tnon_k = dealias*h*Tlinear_k
elif method == 'ETD':
    Tlinear_k = np.exp(h*Loperator_k) 
    def myexp(x):
        if x == 1: return 1.0
        else: return (x-1.0)/np.log(x)
    vmyexp = np.vectorize(myexp) # vectorize myexp (could be jitted)
    Tnon_k = dealias*h*vmyexp(Tlinear_k)
else: print('ERROR: Undefined Integrator')
# Initial condition
rng = np.random.default_rng(12345)
noise = 0.02
n0 = 0.5
n[0] = n0 +noise*rng.standard_normal(n[0].shape)
Noperator_k = n_k.copy() # auxiliary array
nn = n[0].copy() # auxiliary array
n_k[:] = fft2(n[0]) # FT initial condition
# time evolution loop
for i in range(1,Nsteps):
    # calculate the nonlinear operator (with dealising)
    Noperator_k[:] = Noperator_func(nn)
    # updating in time
    n_k[:] = n_k*Tlinear_k + Noperator_k*Tnon_k 
    # IFT to next step
    nn[:] = ifft2(n_k).real 
    # test to output
    if (i % nframes) == 0: n[i//nframes] = nn
\end{lstlisting}

\subsection{2D Phase Field Crystal Model}

\begin{lstlisting}[language=Python]
import numpy as np
from scipy.fft import fft2, ifft2
# PFC model constants
r = -0.25
M = 1.0 # mobility
# Size of the system
N = 2**8 # 2**8 = 256
L = 16*np.pi
x = np.linspace(0,L,N)
dx = x[1]-x[0]
# The time step definition
h = 0.01
T = 1500
Nsteps = int(T/h)
dframes = 1.0 # time step to output
Nframes = int(T/dframes) #frames to the output
nframes = Nsteps//Nframes
# The array of outputs
n = np.empty((Nframes,N,N), dtype=np.float32)
# The Fourier variables
n_k = np.empty((N,N), dtype=np.complex64)
kx = np.fft.fftfreq(N, d=dx)*2*np.pi
k = np.array(np.meshgrid(kx , kx ,indexing ='ij'), dtype=np.float32)
k2 = np.sum(k*k,axis=0, dtype=np.float32)
kmax_dealias = kx.max()*2.0/3.0 # The Nyquist mode
# Dealising matrix
dealias = np.array((np.abs(k[0]) < kmax_dealias )*(np.abs(k[1]) < kmax_dealias ),dtype =bool)
# The linear terms of PDE
Loperator_k = -M*k2*(k2**2-2*k2+1+r)
# The non-linear terms of PDE (with dealising)
def Noperator_func(n):
    return -(k2*M*fft2(n**3))
# Defining the time marching operators arrays
# can be calculated once
if method == 'IMEX':
    Tlinear_k = 1.0/(1.0-h*Loperator_k) 
    Tnon_k = dealias*h/(1.0-h*Loperator_k) 
elif method == 'IF':
    Tlinear_k = np.exp(h*Loperator_k) 
    Tnon_k = dealias*h*Tlinear_k
elif method == 'ETD':
    Tlinear_k = np.exp(h*Loperator_k) 
    def myexp(x):
        if x == 1: return 1.0
        else: return (x-1.0)/np.log(x)
    vmyexp = np.vectorize(myexp) # vectorize myexp (could be jitted)
    Tnon_k = dealias*h*vmyexp(Tlinear_k)
else: print('ERROR: Undefined Integrator')
# Initial condition
rng = np.random.default_rng(12345)
n0 = -0.085
noise = 0.02*np.abs(n0)
n[0] = n0 +noise*rng.standard_normal(n[0].shape)
Noperator_k = n_k.copy() # auxiliary array
nn = n[0].copy() # auxiliary array
n_k[:] = fft2(n[0]) # FT initial condition
# time evolution loop
for i in range(1,Nsteps):
    # calculate the nonlinear operator (with dealising)
    Noperator_k[:] = Noperator_func(nn)
    # updating in time
    n_k[:] = n_k*Tlinear_k + Noperator_k*Tnon_k 
    # IFT to next step
    nn[:] = ifft2(n_k).real 
    # test to output
    if (i % nframes) == 0: n[i//nframes] = nn
\end{lstlisting}

\bibliography{biblio.bib}

\end{document}